\documentclass[]{article} 

\usepackage[
  paper=a4paper,
  left=12.5mm,
  right=25mm,
  top=25mm,
  bottom=50mm,
  bindingoffset=10mm]{geometry}		
 
\usepackage[utf8]{inputenc}         
\usepackage[T1]{fontenc}
\usepackage{lmodern}				
\usepackage{microtype}				
\usepackage{graphicx} 				
\usepackage{url}					

\usepackage{bbm}
\usepackage{amssymb,amsmath,amsthm}
\usepackage[active]{srcltx}
\usepackage{listings}				
\usepackage{mdwlist}				
\usepackage{setspace} 				

\usepackage{todonotes}				
\usepackage{lscape}					
\usepackage{calc}
\usepackage{footnote}				
\usepackage{rotating} 
\usepackage{tablefootnote} 
\usepackage{lineno} \usepackage{array} 
\newcolumntype{x}[1]{!{\centering\arraybackslash\vrule width #1}} 
\usepackage{booktabs} 
\usepackage{longtable} 		
\usepackage{multicol}	
\usepackage{tabularx}
\usepackage{multirow} 
\usepackage{verbatim} 
\usepackage{amssymb} 
\usepackage{eurosym} 
\usepackage{makecell} 
\usepackage{algorithm}
\usepackage{algpseudocode}
\usepackage{longtable} 
\usepackage{xcolor}
\usepackage[normalem]{ulem}
\usepackage{booktabs}
\usepackage{caption}
\usepackage{subcaption}
\usepackage{lscape}
\usepackage{subscript}
\usepackage{placeins} 
\usepackage{color, colortbl}
\usepackage{float}
\usepackage{rotating}

\usepackage[acronym,toc]{glossaries}
\newacronym{TSN}{TSN}{Time Space Network}
\newacronym{MLP}{MLP}{Multi Layer Perceptron}
\newacronym{GA}{GA}{Genetic Algorithm}
\newacronym{ML}{ML}{Machine Learning}
\newacronym{LSTM}{LSTM}{Long Short Term Memory}
\newacronym{LIFO}{LIFO}{Last In First Out}
\newacronym{FIFO}{FIFO}{First In First Out}

\newacronym{VSP}{VSP}{Vehicle Scheduling Problem}
\newacronym{CSP}{CSP}{Crew Scheduling Problem}
\newacronym{VCSP}{VCSP}{Vehicle Crew Scheduling Problem}
\newacronym{CRP}{CRP}{Crew Rostering Problem}
\newacronym{CSRP}{CSRP}{Crew Scheduling Rostering Problem}
\newacronym{MDVSP}{MDVSP}{Multi Depot Vehicle Scheduling Problem}
\newacronym{BPC}{BPC}{Branch-and-Price-and-Cut}
\newacronym{DFS}{DFS}{Depth-First Search}

\newacronym{CBN}{CBN}{Connection-Based Network}
\newacronym{MILP}{MILP}{Mixed Integer Linear Problem}

\usepackage[
    backend=biber,
    style=apa,
  ]{biblatex}

\usepackage{calc}

\usepackage[colorlinks,linkcolor=black,citecolor=black,urlcolor=black]{hyperref} 

\makeglossaries

\begin{document}
	\title{An overview of optimization approaches for scheduling and rostering resources in public transportation}
    \author{Lucas Mertens \and Lena-Antonia Wolbeck \and David R{\"o}{\ss}ler \and Lin Xie \and Natalia Kliewer}
    
    \maketitle
		
	\begin{abstract}
	Public transport is an essential component in satisfying people’s growing need for mobility. Thus, providers are required to organize their services well in order to meet the high demand of service quality at low operational costs. In practice, optimized planning can lead to considerable improvements for providers, customers, and municipalities. The planning process related to public transport consists of various decision problems, of which the providers are usually responsible for vehicle and crew planning. There is a growing body of literature that recognizes the shift from sequential and iterative to integrated solution approaches for these problems. Integrated optimization of several planning phases enables higher degrees of freedom in planning, which allows for operational costs savings and increased service quality. This paper provides an overview of solution approaches for integrated optimization based on operations research techniques for the vehicle scheduling, crew scheduling, and crew rostering problem, extended by a selected number of relevant related approaches from other industries. Therefore, existent optimization approaches are analyzed with regard to different aspects such as mathematical modeling, optimization objective, and method, as well as the source and scope of the data used for evaluation. Additionally, we analyze the problem dimensions that are usually required in practical applications. In doing so, we are able to point out directions for future research, such as a stronger focus on objectives besides cost-minimization like robustness, schedule regularity, or fairness.
    \end{abstract}
	
	\section{Introduction}\label{sec:intro}
	Urbanization in developed and developing countries leads to quickly growing needs for urban mobility. Cities and municipalities face severe challenges in providing the necessary infrastructure to satisfy these needs. Individual motorized traffic is rather part of the problem than of the solution: Traffic jams leading to adverse effects such as long commuting times, frequent accidents, and air pollution are just some of the issues that arise from cities being overflown with cars \parencite{schrank2019urban}. An efficient public mass transportation system can remedy these problems \parencite{ibarra-rojas2015planning}. In addition to the advantage of accommodating the growing mobility demand at lower external effects and costs, public mass transport is safer as well as more resource-efficient than individual transport \parencite{litman2016hidden}.
	
	Traditionally, public transport used to be provided solely by the public sector. However, it has been deregulated in many countries. Nowadays, the transport services offered by private companies are substantial, and competition in this area is increasing \parencite{hrelja2018partnerships}. For a public transport provider, an effective and efficient operation is crucial to face the trade-off between operating costs and service quality. Thus, in each phase of the planning process, this trade-off is considered. Moreover, further goals like schedule robustness, regularity, travel satisfaction, and fairness for employees have gained more attention in recent years \parencite{borndorfer2015duty,cats2016robustness,abenoza2017travel}.
	
	Since the underlying decision problems are not trivial to solve (to optimality), public transport planning has been extensively studied in the literature. Usually, the public transport planning process is divided into planning steps that have to be performed subsequently. However, recent advances in optimization methods allow a gradual integration of the optimization subproblems arising from subsequent planning steps. While better network designing, line planning, and timetabling effect both customer satisfaction and cost structure, vehicle scheduling, crew scheduling, and crew rostering mainly influence the provider’s profit as well as operational timeliness. Superior calculated vehicle and crew schedules lead to lower investments due to less required vehicles and personnel and lower variable costs due to decreased deadheading distances and improved duty allocation. Furthermore, crew rostering impacts costs and employee satisfaction likewise, as crew members desire a fair distribution of duties and workload. Integrating two or three subproblems increases the degree of freedom for these decision problems, and thus, schedule and roster quality may improve. Other industries like railway or aircraft face similar challenges, therefore their solution approaches might be transferable to public transport planning. 
	
	The last decade has witnessed an enormous increase in publications on integrated optimization approaches for public transport planning problems. In 2015, \citeauthor{ibarra-rojas2015planning} conducted a literature review on solution approaches for bus transport systems. In order to extend and update this overview, we will analyze the state-of-the-art approaches that follow different variants of integration and objectives in this paper. We first introduce the operational problems in public transport and point out the contribution of the sequential approach in Section~\ref{sec:process}. Second, we point out the ongoing shift from sequential to an integrated approach in Section \ref{sec:analysis}. 
	
	\section{Decision problems within the operational public transport planning process}\label{sec:process}
	
\begin{figure}
	    \centering
	    \includegraphics[width=0.75\textwidth]{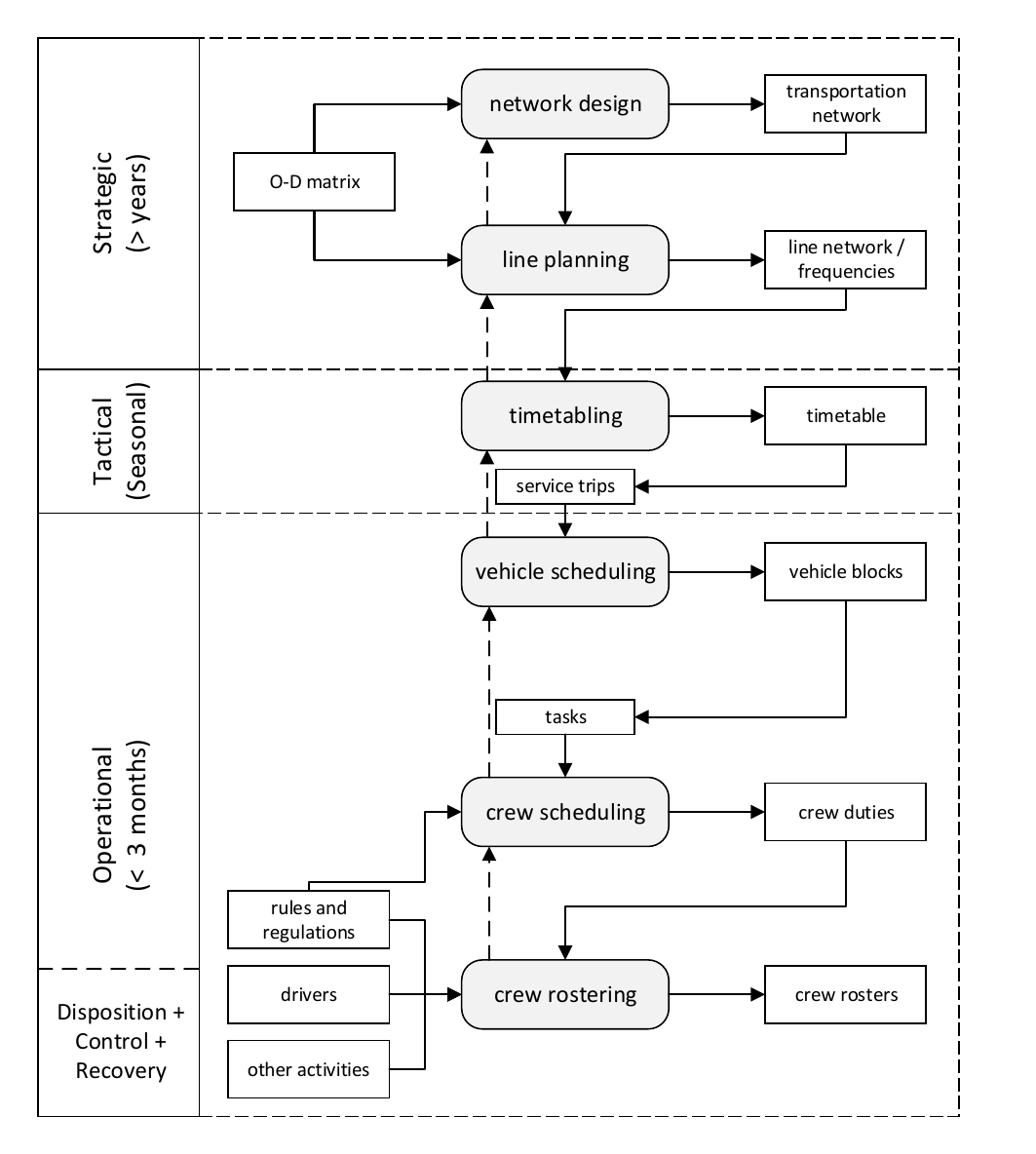}
	    \caption{The sequential planning process in public bus transit, as illustrated in \cite{xie2014decision}.}
	    \label{fig:planprocess}
	\end{figure}

	The planning process in public transport comprises various decision problems, which can be grouped according to their planning horizons (see Figure~\ref{fig:planprocess}). On a strategical level, public transport providers plan long-term, e.g., the network design and the planning of lines (routes, frequencies). Tactical decisions, however, aim to provide timetables and to reduce the operational costs in the medium term \parencite{ibarra-rojas2015planning}. In practice, such decisions are most likely made by the principal, e.g., the municipality \parencite{huisman2004integrated}. We consider strategical and tactical planning decisions regarding routes, frequencies, and timetables as input for the operational planning tasks of vehicle scheduling, crew scheduling, and crew rostering. Therefore, in the scope of this paper, we assume that public transport providers focus on minimizing costs concerning vehicles and staff in the short term when operationally deciding on their transport and employee schedules \parencite{ibarra-rojas2015planning}. Following this, we look at three decision problems: The \gls{VSP}, the \gls{CSP}, and the \gls{CRP}, which are introduced in the following.

	\subsection{Vehicle scheduling problem}
	Given a timetable with specified service trips, the \gls{VSP} relates to generating an optimal vehicle schedule that covers all service trips and achieves the lowest operational costs or optimizes further objectives \parencite{daduna1995vehicle}. A service trip is defined by the line it belongs to, a departure and arrival time, as well as the corresponding locations. To achieve a sequence of compatible trips, additional deadhead trips can be added to connect subsequent service trips. These deadhead trips comprise all unloaded trips, including the departure from (pull-out) and arrival at (pull-in) the depot. A solution to the \gls{VSP} corresponds to a vehicle schedule consisting of vehicle blocks, each representing a feasible sequence of trips for one vehicle \parencite{bodin1981classification}. Thus, a vehicle block comprises one or several vehicle rotations starting from a depot, executing one or more service trips, and returning to a depot.
	
	Solving a \gls{VSP} is not a trivial task and can vary greatly depending on practical requirements and circumstances. The fundamental \gls{VSP} is characterized by a single depot, a homogeneous fleet, and the objective to minimize costs only \parencite{daduna1995vehicle}. One of the first optimal solutions of such a \gls{VSP} originates from \cite{saha1970algorithm}. However, modern \gls{VSP} evolved to cover a more complex environment. Opposed to originating from one depot only, the \gls{MDVSP} considers multiple depots as well as multiple vehicle types and vehicle type groups. This enhancement majorly affects the way of solving the problem. Whereas the single-depot \gls{VSP} is described as a polynomially solvable minimum cost flow problem, the \gls{MDVSP} is considered to be NP-hard \parencite{bertossi1987some}. By utilizing a linear programming approach with column generation, \cite{loebel1999solving} exactly solve the \gls{MDVSP}. Considering multiple depots as well as a heterogeneous fleet, \cite{kliewer2006time} present a \gls{TSN} to efficiently model a network associated to the \gls{MDVSP}. By modeling the \gls{MDVSP} as a \gls{TSN}, the solution space can be reduced significantly. As a result, optimally solving the multicommodity min-cost flow MIP-formulation of the \gls{MDVSP} for real-world instances is made possible. However, not only the underlying problem shifted to a more complex model, but also the objective itself adapted. Whereas in the beginning, the focus was primarily on cost-related objectives, other goals like the schedule robustness are increasingly considered in recent years \parencite{kramkowski2009increasing, naumann2011stochastic}. Different dimensions regarding constraints such as the limited ranges of electric buses \parencite{adler2014routing, reuer2015electric} and objectives shape each \gls{VSP} individually. Several modeling approaches, as well as specialized solution strategies for the \gls{VSP} and its extensions, have been developed in the last decades. For an overview on vehicle scheduling and corresponding solution approaches, we refer to \cite{bunte2009overview} and \cite{pepin2009comparison}. 
	
	\subsection{Crew scheduling problem}
	In sequential planning, the decision problem of crew scheduling arises succeeding to vehicle scheduling. The \gls{CSP} (also known as driver, duty, or shift scheduling) aims at finding a daily cost-optimal duty allocation that encompasses all trips of the vehicle blocks \parencite{borndoerfer2001duty}. These duties are not assigned to specific drivers yet. Each anonymous duty is associated with a predefined generic duty type. These heterogeneous duty types are characterized by different lengths and attributes. A duty type, e.g., considers legal requirements on working and break times as well as company-specific regulations such as the kind of qualifications required \parencite{freling2004decision}.
	
	Due to the vast amount of possible solutions based on the predefined duty types covering the vehicle schedule's trips, solving the \gls{CSP} is considered to be NP-hard \parencite{fischetti1987fixed}. The complexity of finding a solution to the \gls{CSP} correlates with the number of trips and especially with the quantity and diversity of the generic duty types. Depending on practical requirements, each duty type at least considers legal, union-related, and company-defined regulations. These characteristics vary significantly regarding each problem specification. In solving the \gls{CSP}, it has prevailed to split all vehicle blocks into segments according to predefined relief points \parencite{desaulniers2007public}. Relief points indicate locations at specific times, which allow an exchange of drivers. The tasks between two relief points represent the smallest unit of work that has to be covered by the same driver and is called duty element. Combining consecutive duty elements and adding sign-on and sign-off tasks results in a possible shift, and is called a piece of work. Final duties are composed of one or more pieces of work, where usually two pieces of work are separated by a break \parencite{desaulniers2007public}.
	
	Since the emerging duties are not associated with specific drivers yet, commonly cost criteria shape the objective of solving the \gls{CSP} \parencite{ernst2004staff, huisman2004integrated}. Depending on the practical application, both minimizing the total amount of daily duties as well as minimizing the total required work time are achievable tasks. Whereby the former determines the minimum demand of employees on a daily basis, the latter aims at an optimal duty structure by avoiding unnecessary breaks or waiting times. Utilizing fixed costs for duties and an hourly rate for the working time, these objectives are usually transformed into one that minimizes the total costs \parencite{desaulniers2007public}. Commonly, for solving the \gls{CSP} a column generation approach in combination with Lagrangian or LP-relaxation considering a set covering or partitioning problem is utilized. Solution approaches for crew scheduling are reviewed in detail in \cite{huisman2004integrated} and \cite{ernst2004staff}.
	
	\subsection{Crew rostering problem}
	
	The crew rostering (or driver rostering) is concerned with assigning anonymized duties to specific drivers. The results are individual schedules for every crew member, so-called crew rosters \parencite{freling2004decision}. As opposed to crew scheduling, where the foremost objective is to minimize operative costs, crew rostering takes crew welfare, such as balancing workload and additional individual characteristics of each crew member and efficiency objectives, e.g., minimizing layovers and crew deadheading, into account. Complementing the  legal daily duty requirements, already respected within the \gls{CSP}, further law and labor union rules have to be considered, solving a \gls{CRP}. These additional requirements range from minimal break times between two consecutive shifts to a maximum weekly workload for a single driver. In constructing personalized schedules, two different kinds of crew rosters can be distinguished, namely cyclic and non-cyclic rosters \parencite{xie2015cyclic}. The cyclic roster occurs to be the less sophisticated approach and is developed for a group of drivers with similar qualifications and preferences. A regular, repeating working pattern is established for the entirety of drivers. This pattern is constructed as such that all legal requirements are met, and the workload is allocated evenly. However, this roster occurs not to feature a high degree of individuality. Non-cyclic rosters, on the other hand, offer the possibility to develop personalized schedules for a medium to long period of time \parencite{xie2015cyclic}. Depending on the extent to consider individual preferences and shift requests, constructing a non-cyclic pattern requires sophisticated techniques. A multicommodity network flow formulation is developed in \cite{xie2015cyclic} to deal with both cyclic and non-cyclic rostering, also in \cite{mesquita2015decompose} for non-cyclic rostering. In order to deal with the complexity, (Meta-)heuristics are applied to solve non-cyclic rostering, such as in \cite{xie2017metaheuristics}, \cite{mesquita2015decompose}. \cite{ernst2004staff} and \cite{vandenbergh2013personnel} cover the crew rostering problem in their literature reviews and elaborate approaches to solve the \gls{CRP} utilizing both cyclic and non-cyclic rosters. As a first step towards more robustness in crew rostering, \cite{xie2012integrated} consider a simplified version of rostering but incorporate possible reserve shifts to cover the absences of drivers.

	\section{Partial Integration \& Integrated approaches} \label{sec:analysis}
	
	
	The three decision problems -- more precisely \gls{VSP}, \gls{CSP} and \gls{CRP} -- have been extensively studied by scholars. Various methods have been proposed to find optimal or close-to-optimal solutions to each of these problems \parencite{bunte2009overview, ernst2004staff, vandenbergh2013personnel}. These problems constitute consecutive phases \parencite{desaulniers2007public} within the operational public transport planning process. Thus, choosing a sequential approach to solving the entirety of these problems is straightforward. In such an approach, the output of the previous phase is used as an input for the subsequent planning problem. However, this traditionally utilized approach may not lead to a globally optimal solution. A slightly adjusted timetable, e.g., might lead to more freedom for solving the \gls{VSP} and hence a lower demand for buses. Here, the gain from consecutive steps can outweigh the loss of the adjusted prior phase or, due to choosing an indifferent solution of a previous step, a Pareto-efficient improvement might even be possible. As a result, iteratively solving the three sequential phases in order to leverage knowledge gained in every iteration can improve the overall solution. However, repeated executions of each phase might lead to prohibitively long run times or, due to a fixed number of iterations, to local optima. Opposed to sequential or iterative approaches, which solve each of the problems separately, integrated approaches solve the \gls{VSP}, \gls{CSP}, or \gls{CRP} conjointly. As a result, superior solutions are attainable within acceptable computation time, even for problem instances of realistic size. 
	
	We distinguish between the integration of the first two phases (\gls{VSP} + \gls{CSP} in the following referred to as VCSP) and the last two phases (CSP + \gls{CRP} in the following referred to as CSRP). The highest level of integration is achieved by simultaneously considering all three decision problems (\gls{VSP} + \gls{CSP} + \gls{CRP} in the following referred to as VCSRP).
	
	The number of publications for integrated solution approaches is unevenly distributed. In contrast to the wide range of publications considering the VCSP, there are only three approaches for the VCSRP. The main objective in either integrated approach is usually minimizing costs -- while in recent years, additional objectives such as robustness, regularity, and fairness have become increasingly important. Similar to publications covering the VSP only, there exists an evident trend to modeling the underlying problem as a \gls{TSN} instead of a connection-based Network (CBN). Due to the integration of the planning phases, many approaches use column generation and (meta-)heuristics (such as genetic algorithms, simulated annealing, ant colony algorithms) to solve the remaining complexity problem. More than two-thirds of the evaluated solution approaches use real data for evaluation and thus examine the applicability of the methods in practice. 
	It is noteworthy that the majority of approaches employ combinations of solution approaches instead of individual exact or heuristic methods. Looking at the methods, special attention is paid to column generation, as it is prevalent in the sample, as well as non-exact heuristics that are used. The right choice of model and combination of solution algorithms facilitates solving problem instances of realistic size. However, within the regarded sample, only a few publications from the bus industry solve VSP instances of realistic urban size \parencite[e.g.,][]{kliewer2012multiple, amberg2018robust, steinzen2010time}. Many rely on evaluation using the random benchmark instances published in \cite{huisman2004integrated}.
	
	Similar decision problems under consideration occur in several industries. Three are identified as the major industries: Airline, railway, and public buses. Regarding vehicle scheduling, similarities, as well as differences, are evident. All three industries have the goal in common to minimize the operational costs and utilize the least possible number of vehicles. However, the details of either industry differ greatly. Due to high initial costs for rolling stock and railroads, as well as long construction times for the latter, planning in the railway industry is highly constrained by its infrastructure. In contrast, vehicle scheduling for a public bus provider offers more decision-making possibilities. Various existing roads can be used, and different vehicle types offer higher degrees of freedom in planning. Given a fixed number of vehicles, scheduling for the airline industry is the least restricted one. Changing the route of an airplane is usually only restricted by costs, but not by infrastructural conditions. Depending on the preconditions of each unique industry, mathematical modeling might be more challenging given increased infrastructural requirements. The number of constraints correlates strongly with the model's degrees of freedom. More flexibility in planning leads to an increased solution space. Both the quantity of constraints and the size of the solution space enable different solution approaches and might lead to different expedient ways of solving the specific planning problem. Similar to vehicle scheduling, both crew scheduling and crew rostering share similarities but differ in detail. As previously described, labor law and other legal provisions, as well as collective and individual agreements, restrict the \gls{CSP} and \gls{CRP} within the mentioned industries \parencite{freling2004decision, guo2006partially, xie2015cyclic}. However, buses, e.g., only need one driver while airplanes and trains must have a crew. Crews typically consist of more than two members who have to fulfill specific tasks and functions, and are thus typically planned as teams \parencite{freling2004decision}. Compared to public bus transport, the railway and airline industry possibly cover huge distances. Thus, the crew rostering has to consider individual home bases, take lodging into account and return each crew member to its origin at some point \parencite{wen2021airline}.
	
	Most studies in our sample deal with the public bus transport industry (in either urban or rural environments). There are some important exceptions from other industries, especially concerning \gls{CSRP} such as \cite{sandhu2007integrated}, \cite{freling2004decision}, \cite{guo2006partially}, \cite{medard2007airline} and \cite{souai2009genetic} in the airline industry and \cite{freling2004decision} as well as \cite{borndoerfer2012rapid} in railway.
	
	In the following sections, we discuss the solution approaches from the literature concerning the pairwise integrated problems (i.e., the VCSP and the CSRP), and the “fully” integrated problem (i.e., the VCSRP) in more detail. 
	
	\section{Pairwise integrated optimization}
	\subsection{Integrated vehicle and crew scheduling}\label{sec:VCSP}
	The majority of solution approaches for the VCSP in our sample follow a column generation scheme to generate vehicle schedules and anonymous duties for a given timetable and corresponding service trips. The VCSP is the master problem, and duties are generated as columns by solving the pricing problem as a constrained shortest path problem. All approaches investigated have in common that minimizing costs is the central objective criterion. In recent years, further optimization objectives such as robustness \parencite{huisman2004robust, amberg2011increasing, kliewer2012multiple, amberg2018robust} and schedule regularity \parencite{steinzen2009branching, amberg2011approaches, amberg2012robuste} have been considered. 
	
	For the corresponding \gls{VSP}, the underlying network is usually explicitly modeled. Historically, integrated optimization approaches have focused on using a CBN with depots and stops as nodes, and all possible connections, including pull-ins and pull-outs, are enumerated as arcs, such as in \cite{ball1983matching}, \cite{freling1999overview}, \cite{friberg1999exact}, \cite{gaffi1999integrated}, \cite{freling2001applying} and \cite{freling2003models}. This approach might be most intuitive and was used mostly in the last century. Recent network modeling approaches shift towards a \gls{TSN} where time-space nodes represent possible arrivals and departures at a location and where only feasible connections are modeled as arcs such as in \cite{gintner2008crew}, \cite{keri2008simultaneous}, \cite{steinzen2010time}, \cite{amberg2011approaches}, \cite{amberg2011increasing}, \cite{kliewer2012multiple} and \cite{amberg2018robust}. The \gls{TSN} method has the advantage that much fewer connections are included, which reduces the model complexity tremendously -- especially for larger instances. \cite{gintner2005solving} report that the number of arcs in the \gls{TSN} amounts to 1-3\% of all arcs in an equivalent CBN. Thus, the problem size could be reduced significantly without reducing the solution space because all compatible trips are implicitly connected.
	
	\subsection{Integrated crew scheduling and rostering}\label{sec:CSRP}
	In the airline and railway industries, crew scheduling and crew rostering are usually considered sequentially (see \cite{caprara1999solution, lee2003scheduling, yunes2005hybrid, medard2007airline} since it is not yet possible to find an optimal solution for one of the two planning steps with current state-of-the-art technologies for realistically sized models. An overview of the developments until 1998 for air and rail transport was presented in \cite{ernst2001rail}.
	
	Integrated planning has received increasing attention since the 2000s, with a focus on airlines and railways. Due to the high combinatorial complexity of integrated planning, approaches to partial or iterative integration were first published. In \cite{ernst2001integrated} the number of paired personnel crews in crew scheduling is taken into account. 
	
	Most integrated crew scheduling and crew rostering approaches deal with airline optimization \parencite{freling2004decision, guo2006partially, medard2007airline, souai2009genetic} and only a few tackle the public bus transit \parencite[e.g.][]{xie2012integrated, xie2013column, xie2017metaheuristics, xie2015cyclic} and the railway industry \parencite[e.g.][]{borndoerfer2014integrierte, lin2019integrated}.
	
	An iterative method through a feedback mechanism between the \gls{CSP} and \gls{CRP} is implemented in \cite{caprara2001global}. All duties are generated in the first phase, and the number of duties is reduced by heuristics in the second phase, such that instances with various compositions with real-world characteristics can be solved. \cite{guo2006partially} focus on partial integration based on the aggregated \gls{TSN}. In the first step, instead of a single duty, a chain of duties is generated, taking into account the individual activities of the crew members planned in advance. In addition, the number of crews is also taken into account in this step \parencite{guo2006partially}. The approach can solve even instances of up to 1977 tasks, considering 188 crew members, in acceptable time ($\sim$ 15.5 minutes).
	
	In \cite{zeghal2006modeling} the integration problem is formulated as an integer linear program, and a new heuristic method is developed, which is used in a search procedure for a subtree based on a rounding strategy. A decision support system is developed in \cite{freling2004decision} for integrated crew scheduling in the airline and railway sector, and a general set partitioning model is formulated, and a state-of-the-art branch and price solver is generated. In \cite{saddoune2011integrated} and \cite{saddoune2012integrated} a column generation approach is used to reduce the computing time of the sub-problem.
	
	In further research projects regarding the complete integration of the two planning phases, meta-heuristics, in particular specialized genetic algorithms, are successfully used to solve the integrated problem \parencite[see][]{souai2009genetic, chen2013integrated}.
	
	Because of lower operational costs, optimization approaches for public transport were developed only about ten years later than in the integrated planning for the airline and railway industry. A Bender's decomposition approach is used in \cite{borndoerfer2014integrierte}, where the crew rostering was simplified in such a way that the duty sequences are anonymous and shift and duty templates were used instead of services. In \cite{xie2015decision} it is shown that in practice, it is often critical to underlay shifts with concrete duties.
	
	\cite{lin2019integrated} propose a \gls{BPC} algorithm for solving the \gls{CSRP} for the Taiwanese railway system with regard to standby personnel. They compare the results with solution approaches using expert knowledge or rules of thumb, commercial standard solvers for the associated \gls{MILP} and a sequential \gls{DFS} based algorithm for several instances reaching real-world problem sizes regarding the number of tasks to be performed. The employed \gls{DFS} first enumerates all potential duties, then identifies the minimum required duties to cover all tasks as a set partitioning problem, and finally solves the shift-assignment to optimality. Only the \gls{BPC} algorithm was capable of solving all instances, whereas Gurobi and the \gls{DFS}-based algorithm are only tractable for the smallest and second-smallest instance, respectively.
	For the smallest instances, the \gls{BPC} can recreate the optimal solution in less time and is the only algorithm capable of solving all problem instances.
	
	In addition to cost minimization, younger approaches aim at optimizing for further goals such as the maximization of fairness of the drivers' shift allocation and the regularity of duty rosters to increase satisfaction \parencite[e.g.][]{borndoerfer2017integration, quesnel2020improving}).

	\section{Integrated vehicle and crew scheduling and rostering}\label{sec:VCSRP}
	
	Few publications look into integrating all three phases, all of which take up the bus industry.
	\cite{shen2009integrated} consider several data sets and circumstances. They use data from the Beijing Bus group to point out the practical constraints that derive from Chinese law and culture. These include built-in meal periods, multi-type bus scheduling, and restricting drivers to one or two particular buses. The authors develop an iterative sequential heuristic algorithm that consists of three steps: Firstly, the \gls{VSP} is solved with a local search based on $n$-opt operators. Then, the \gls{CSP} is solved using a tabu-search heuristic. Finally, driver rosters are proposed to the user and can be modified through an interface. According to the authors, it is possible to find feasible solutions for instances up to 107 buses and 164 duties within an appropriate time frame of some minutes. The authors report savings in vehicle costs close to 4.5\% and driver wages of approximately 9.9\% when comparing with manually built solutions.
	
	\cite{mesquita2011new} use data from a bus company in Lisbon to demonstrate their preemptive goal programming-based heuristic approach that prioritizes the \gls{VCSP} over the \gls{CRP} part. Their approach is able to generate optimal solutions within a short computing time for most instances. When considering all costs, however, some instances could not be solved within a reasonable time limit. Their integer formulation consists of a preemptive goal programming framework that prioritizes the integrated vehicle-crew-scheduling goals over the driver rostering goals. The problem is first decomposed to solve one \gls{VSP} + \gls{CSP} per day and then establish a roster for a longer time horizon.
	
	Two years later, \cite{mesquita2013decomposition} manage to outperform the traditional sequential approach by integrating \gls{VSP}, \gls{CSP} and \gls{CRP} with a Bender's decomposition problem formulation using a multicommodity network flow formulation, set covering and covering-assignment elements. They tackle the integrated problem by dividing it into a master problem that contains the \gls{VSP} and \gls{CSP} and a sub-problem for the \gls{CRP}. Information from the sub-problem and its dual solution are used to find better duties for the \gls{CSP}. The authors minimize vehicle and driver costs and take into account constraints regarding roster balancing and coverage of all daily duties. Using data from two bus companies in Portugal, their rosters have to match predefined days-off patterns based on the requirements of these companies. The planning horizon for rosters is seven weeks long.
	
	All three papers, \cite{shen2009integrated}, \cite{mesquita2011new}, and \cite{mesquita2013decomposition}, evaluate their proposed algorithms on real-world instances. However, they are too small (108 to 238 timetable trips) to represent a larger, realistic urban bus system.
	\\
	In summary, public transport providers recognize the necessity to organize their services efficiently. Because of increasing urbanization, the demand for public transport is rising, and thus competition and the need for efficiency in the public transport planning process rise as well. Integrating the operational phases of \gls{VSP}, \gls{CSP}, and \gls{CRP} gives public transport providers more degrees of freedom and can lead to better schedules and rosters.
	
	Our findings include that there are three forms of integrated problems that are solved in the literature. While most authors focus on the \gls{VCSP}, some approaches consider the \gls{CSRP} and a few recent publications tackle the challenge of solving the VCSRP.
	
	All approaches to solving the VCSRP deal with the bus industry. This might be because crew rostering is much easier when considering only one driver, as compared to multiple-person crews in airline or railway planning, which would exacerbate the integration even more.
	
	Moreover, the majority of scholars focus on minimizing costs in their approaches. However, a few authors have considered other objectives such as robustness (see an overview of different robustness approaches in public transit in \cite{ge2022robustness}), regularity, and fairness of schedules in recent years. This indicates that diverse objective functions are becoming more common over time. In a recent study, \cite{ge2022revisiting} showed that adding a robustness objective to the VCSRP model of \cite{mesquita2013decomposition} does not take much additional computational time in this application, which is an interesting result. Many standard combinatorial optimization models are used for decision problems, including the minimum cost flow problem and its various special cases (e.g., the resource-constrained shortest path problem, the multicommodity flow problem, and the linear assignment problem). Set partitioning formulations are often used to solve the \gls{CSP}. Some authors prefer the easier set covering formulation where crew members become passengers in the case of overlapping assignments. The \gls{TSN} formulation is a powerful tool to reduce network size and thus computing time compared to \gls{CBN}, where all deadhead trips are explicitly modeled.
	
	In terms of solution techniques, column generation stands out as the most powerful operations research method to solve integrated decision problems. It is usually accompanied by relaxation techniques. Lagrangian relaxation seems to work best for quickly finding reasonable bounds for the integer solution. Branch-and-bound and branch-and-price techniques are popular to find feasible integer solutions. Heuristics and meta-heuristics are used to speed up the solution process. They include tabu-search, simulated annealing, ant colony algorithms, genetic algorithms, and smaller-scale greedy heuristics.  
	In conclusion, we can say that there exists no single best approach to integrated solve the \gls{VSP}, \gls{CSP}, and \gls{CRP}. Which solution approach yields the best results is always subject to the specific problem settings. Among others, it is important to take into account the source and nature of data that is used, the size of the instances, and the relevant constraints. Every new approach can be a game-changer for some situations, while in others, it might prove less useful.
	
	\section{Future Research}\label{sec:future-research}
	
	For future research, we propose to focus on further integrating the public transport planning process. Since there are only three approaches towards a threefold integration of all operational phases and the first results are promising, more effort is needed in this direction. In addition, strategic decision problems may also be included in integrated planning. In the literature, the first approaches towards integrating timetabling or vehicle routing with vehicle scheduling can be found. The long-term aim of an integrated public transport planning process where all sub-problems are solved simultaneously and thus making it possible to use all degrees of freedom is still a long way off. 
	
	At the same time, the existing approaches should be enhanced in order to make them suitable for real-world use with larger instances and more complex data sets, such as public transport providers in larger cities. Furthermore, schedule robustness, regularity, crew or driver preferences, and fairness as optimization objectives should be examined more closely. There are some first steps taken in individual publications in integrated scheduling. The crew scheduling literature offers many more starting points that could be incorporated into integrated planning as well. More publications that explicitly deal with these topics are desirable.
		
	\printbibliography[heading=subbibliography]

\end{document}